\documentclass[runningheads]{llncs}
\usepackage{graphicx,bm}
\usepackage{amssymb}
\usepackage{algorithm,algpseudocode}
\usepackage{caption}
\usepackage{subcaption}
\begin{document}
\def\eh{\hat{e}}
\def\P{\mathbb{P}}
\title{BOP-Elites, a Bayesian Optimisation Algorithm for Quality-Diversity Search 
}
%
%
\author{Paul Kent\inst{1,2} \and
Juergen Branke\inst{3,4}}
\authorrunning{P. Kent \& J. Branke}
%
\institute{Mathematics of Real World Systems, University of Warwick, UK \and
\email{Paul.Kent@warwick.ac.uk}\\
\url{https://warwick.ac.uk/fac/sci/mathsys/people/students/2018intake/kent/} \and
Warwick Business School, University of Warwick, UK \and
\email{Juergen.Branke@wbs.ac.uk} \\
\url{https://www.wbs.ac.uk/about/person/juergen-branke/}}

\maketitle              
\begin{abstract}
Quality Diversity (QD) algorithms such as MAP-Elites are a class of optimisation techniques that attempt to find a set of high-performing points from an objective function while enforcing behavioural diversity of the points over one or more interpretable, user chosen, feature functions.

In this paper we propose the Bayesian Optimisation of Elites (BOP-Elites) algorithm that uses techniques from Bayesian Optimisation to explicitly model both quality and diversity with Gaussian Processes. By considering user defined regions of the feature space as 'niches' our task is to find the optimal solution in each niche. We propose a novel acquisition function to intelligently choose new points that provide the highest expected improvement to the ensemble problem of identifying the best solution in every niche. In this way each function evaluation enriches our modelling and provides insight to the whole problem, naturally balancing exploration and exploitation of the search space. The resulting algorithm is very effective in identifying the parts of the search space that belong to a niche in feature space, and finding the optimal solution in each niche. It is also significantly more sample efficient than simpler benchmark approaches. BOP-Elites goes further than existing QD algorithms by quantifying the uncertainty around our predictions and offering additional illumination of the search space through surrogate models. 

\end{abstract}
\section{Introduction}
In many optimisation problems, decision makers are looking for a variety of solutions. In the problem setting we consider, besides the objective function to be optimised, every solution is characterised by a set of feature values. The decision maker a priori categorises these feature values into classes (which we call 'niches' in this paper following the MAP-Elites literature), and is looking for the best solution in each niche.
A solution's feature values are computationally expensive to compute or are linked to the objective function evaluation of the solution. For example, a solution may be evaluated by running a simulation model, and the output of the simulation model is not only the objective function value, but also the feature values of the solution. An application could be in robotics, where we try to identify motion control parameters that will get the robot as quickly as possible to each of a set of target regions. We are then looking, for each of the target regions (niches), for the control parameters that move the robot to the target region in the least amount of time (objective function).

The problem is similar, but different to other classes of optimisation problems where evolutionary algorithms have been successful:
\begin{itemize}
    \item In multimodal optimisation, one is looking for several local optima in the search space, but there is only the objective function (no features that define separate niches, the niches are implicitly defined by the location of local optima).
    \item In multi-objective optimisation one is looking for several solutions with different trade-offs in objectives, objectives have to be optimised simultaneously, but there are no pre-defined niches. 
    \item In multi-task optimisation, several related optimisation problems have to be solved simultaneously. Our problem setting could be framed as multi-task optimisation, where each task is to find the best solution belonging to a particular niche.
\end{itemize}

The problem of identifying a diverse set of solutions according to some additional features is known as quality diversity search, and a popular algorithm is MAP-Elites \cite{cully2015robots} which is based on evolutionary computation.

In this paper, we propose a Bayesian optimisation algorithm for this problem. Bayesian optimisation is a surrogate-based black-box optimisation technique that has recently gained a lot of attention in machine learning, in particular for hyperparameter tuning. The algorithm we propose, called 'BOP-Elites', simultaneously builds surrogate models for the objective function as well as for the feature value functions, and then iteratively decides what solution to evaluate next in order to maximise the expected improvement in final solution set quality. We demonstrate the power of the proposed algorithm on a simple benchmark problem and compare it with simpler versions that consider the problem as a multi-task optimisation problem and either allocate computational budget in a round-robin fashion to the different tasks, or solve each task independently.





The paper is structured as follows. We start with a summary of related work in Section~\ref{sec:related}, followed by a more formal problem definition in Section~\ref{sec:probdef}. Section~\ref{sec:bayesopt} provides a brief introduction to Bayesian optimisation, and the new BOP-Elites algorithm is introduced in Section~\ref{sec:bopelites}. Empirical results are reported in Section~\ref{sec:results}, and the paper concludes with a summary and some suggestions for future work.

\section{Related work\label{sec:related}}
\subsection {Quality-Diversity Algorithms}

Whilst standard optimisation focuses on finding the best solution over an input domain, Quality-Diversity (QD) Algorithms such as Novelty Search (NS) \cite{lehman2008exp} and MAP-Elites \cite{mouret2015illuminating} are approaches that attempt to create a diverse set of solutions and reflect the nuanced way in which natural evolution diversifies as it optimises. 

Like Multi-Objective Optimisation (MOO), NS returns a set of high performing points. While MOO tries to find the Pareto Front, the trade off between two or more user chosen objectives, NS goes against the wisdom of conventional optimisation by focussing on behavioural diversity without any pressure to improve the objective output. Surprisingly, in some domains NS can sometimes find the optimum even without explicitly searching for it. In particular, QD algorithms often perform well in optimising functions that have some level of deception. QD algorithms have even been implemented in search problems that have otherwise been considered too difficult \cite{cully2013behavioral,lehman2011evolving}. 

In MAP-Elites, an evolutionary computation-based approach to QD \cite{cully2015robots}, behavioural classifications are pre-defined instead of being emergent throughout the search process. A behaviour space is separated into discrete niches and whenever a new point is evaluated it is assigned to one of these niches. The best performing solution that ends up in that niche is considered as an 'elite'. Through a process of evolution where points breed with those from other niches, and niches to breed with are chosen uniformly at random, new points are generated that are naturally distributed around the search space. 

There are many real world situations where producing a diverse set of high performing solutions may be desirable. Clear examples exist in robotics, the original field of application for MAP-Elites, but extend far beyond \cite{pugh2016quality}. Applications have emerged for instance in the games industry \cite{fontaine2019mapping,ecoffet2019go} and in computer science where it was used to produce a set of images to trick deep neural networks \cite{nguyen2015innovation}. With development of the field, the range of applications for such algorithms is likely to increase.

\section{Problem definition\label{sec:probdef}}

The BOP-Elites algorithm is designed to  work over a constrained, $d$ dimensional, search domain $X \subset {\rm I\!R}^d$. Each dimension is constrained between a predefined lower and upper limit $[L,U]$. 

We define an objective function over the search domain which returns a scalar objective valuation $f(x) = y$.

\[f(x): X   \to  {\rm I\!R} \]  

We are also given a number of feature functions over the search domain  $X$ which define an $m$ dimensional feature space.

\[{\rm \!G(x)} = \{g_1(x),..., g_m(x)\} : X   \to  {\rm I\!R}^m \] 

Each feature function returns a scalar value $g_i(x)=z\in {\rm I\!R} $ which describes some high level characteristic. The idea is that these functions model some interpretable feature that can be classified into various niches, representing the behavioural diversity of the search space. Typically these features are also time consuming to compute and are linked to the evaluation of a solution. 

For each feature $g_i$, human decision makers specify a priori a set of $b_i-1$ boundaries  defining $b_i$ regions within the dimension of the feature $i$. As feature space dimensions are orthogonal, the boundaries divide the search space into $C = \prod_{i=1}^m b_i$ niches, areas with distinct behavioural qualities. A simple example would be having 2 feature functions each with one boundary at a value of 2. $g_1$ has two region labels 'Slow' ($g_1(x)<2$ ) and 'Fast' ($g_1(x)\ge 2$), and $g_2$ has two labels 'Weak' ($g_2(x)<2$) and 'Strong' ($g_2(x)\geq 2$). This leads to $C=4$ niches with the interpretable labels $c_1$ 'Slow and Weak', $c_2$ 'Fast and Weak', $c_3$ 'Slow and Strong' or $c_4$ 'Fast and Strong'. 

The goal of the BOP-Elites algorithm is to find  the set $S^*$ of optimal solutions within each of the predefined niches in as sample efficient a way as possible. Formally this can be expressed as
\[S^* = \{e_1^*,...,e_{C}^*\}\]
\[e_i^*= argmax_x\left( f(x) \right) \,\,\,\, s.t \,\,\,\, x \in X,\,\, {\rm \!G}(x) \to c_i  \]

A solution provided by any algorithm will consist of a set of 'elites' $\hat{S} = \{\hat{e}_1,...,\hat{e}_{C}\}$, the best solution found for each niche. If the algorithm couldn't identify any point in a niche $c$, then $\hat{e}_c$ may be empty.

Suggested solution sets are evaluated by the following total error function, calculated as the distance from the true ensemble optima
\begin{equation}
     TE=\sum_{i=1}^{C}\left( f(e^*_i)-f(\hat{e}_i) \right) \label{TE}
\end{equation} 
 If $\hat{e}_i$ is empty some low default fitness is assigned to niches where no solution has been found, i.e. $f(\hat{e}_i)=f_{min}$.

\section{Bayesian Optimisation\label{sec:bayesopt}}

Bayesian Optimisation (BO) is a derivative free strategy for optimisation of expensive black-box functions. Since the objective function is expensive, BO builds a cheap-to-evaluate surrogate model and performs an extensive search on this model before selecting a new point to sample from the true function. Typically, BO treats the objective function as coming from a probability distribution and places a prior distribution over function space, using sequential samples from the objective function to update the prior into a posterior distribution.  

Key to BO algorithms is the acquisition function (sometimes called infill criterion), a calculation performed on the posterior distribution, which attempts to predict the value of sampling new points in the input space. Acquisition functions are problem specific and designed to balance exploration and exploitation in the search. Examples of acquisition functions include probability of improvement \cite{Kushner64}, Expected Improvement (EI) \cite{Jones98} and Knowledge Gradient (KG)\cite{Scott11}. 

The most widely used BO surrogate model, and the one implemented in BOP-Elites, is the Gaussian Process (GP) model \cite{RassmussenWilliams}. Assuming that the latent function $f(.)$ can be suitably modelled with a GP, for any finite set of observed data pairs $\tau^n= \{(x_1,y_1),...,(x_n,y_n)\}$, $X_n = \{x_1,...,x_n\}, Y_n = \{y_1,...,y_n\}$ we model this function as a multivariate Gaussian random variable. Assuming a prior mean $\mu^0= 0 $ we predict a new function call $f(x)$ from an unobserved point with posterior mean and variance 

\begin{equation}
    {\rm I\!E} \left[ f(x)| \tau^n \right]= \bar{f}(x) = K_*K^{-1}Y \label{GP1}
\end{equation} 
\begin{equation}
    var(f(x)) = s(x)^2 = K^{T}_{**}-K_*K^{-1}K_* \label{GP2}
\end{equation} 

where $K$ is the kernel, or covariance, function. $K_* = K(x,X_n)$, $K^{-1} = K(X_n,X_n)^{-1}$ and $K_{**}=K(x,x)$. While many kernels exist, throughout this paper we use the standard RBF kernel in the BOP-Elites algorithm. This derivation is for the modelling of deterministic functions and therefore does not include a noise term, but an adaptation to noisy functions would be straightforward. GP regression is covered in greater depth in the influential work of Rassmussen and Williams \cite{RassmussenWilliams}.  

\section{BOP-Elites Algorithm\label{sec:bopelites}}

The BOP-Elites algorithm is a model-based optimisation algorithm using Gaussian Process (GP) surrogate models. Both the objective and feature functions are modelled with GPs and used to infer the values and behavioural characteristics of points in the search domain. At each iteration an acquisition function is used to select the next point to be sampled.

Once evaluated, the new observation is used to enhance our surrogate models and the whole process proceeds until the search budget has been exhausted. This approach is particularly suited to the optimisation of expensive black-box functions and is designed to be  sample efficient.

\begin{figure}
\includegraphics[width=\textwidth]{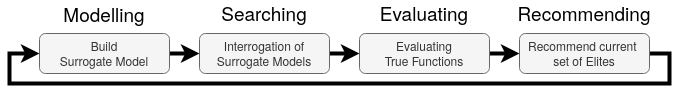}
\caption{The flow of the BOPElites Algorithm.} \label{fig1}
\end{figure}


In our  experiments, BOP-Elites works over a discretised set of points  $X_d \subset{X}$. This choice was made for simplicity and reproducibility but the acquisition function is continuous and future implementations can perform true optimisation over the continuous input domain. It should be noted that, in this instance, BOP-Elites solves the approximate problem of finding the best solutions within the discretised domain. 

The algorithm is initialised with a set of $n_0$ initial points sampled uniformly at random from $X_d$. The objective and feature functions are evaluated for each initial point and stored as a list of observations, the data $D=\{x_i,y_i,G(x_i)\}_{i=1}^n$. We construct the list of 'elites' $\hat{S} = \{\eh_1,..,\eh_{C}\}$, the highest performing solution observed so far for each niche, which is updated whenever a point is sampled that outperforms the elite in its respective niche. Where there is currently no elite for a niche, meaning no points have been observed in this niche, BOP Elites considers its objective contribution as $f_{min}=0$. The list of elites is then the set of points recommended by the algorithm at any iteration. It is also fundamental to the calculation of the acquisition function which is discussed next.

\subsection{BOP-Elites Acquisition Function\label{acqfun}}

The goal of our algorithm is to find the optimal solution for each niche. This requires the simultaneous solving of two problems:
\begin{enumerate}
        \item Given a solution $x$, predict which niche $c$ solution $x$ belongs to.
    \item Given a solution $x$ estimate the improvement gained from $f(x)$ over the current elite in $c$, $f(\eh_c)$.
\end{enumerate}
The second problem can be tackled by modifying the well established EI algorithm \cite{Jones98} designed for global optimisation, here re-purposed to estimate the positive improvement of a new point $x$ over the current elite.

\begin{eqnarray}
    EI_c(x) &=& \mathbb{E}\left[max(f(x)-f(\hat{e}),0)\right] \\
 &=& (\bar{f}(x)-f(\hat{e}))\Phi\left(\frac{\bar{f}(x)-f(\hat{e})}{s(x)}\right)+s(x)\phi\left(\frac{\bar{f}(x)-f(\hat{e})}{s(x)}\right)  \label{nicheimp}
\end{eqnarray}

where $\Phi$ and $\phi$ are the standard cumulative normal  and density functions, $\bar{f}(x)$ and $s(x)$ are the posterior mean prediction and standard deviation at $x$ from the objective surrogate model. EI suggests points that have a high probability of yielding improvement either because the surrogate model predicts an improvement with high certainty or because a point has high uncertainty that could yield a high positive value. 

Implementing EI requires us to choose which elite to compare to. It is entirely possible for the elite in two neighbouring niches to have very different objective values and without allocating the solution candidate to the correct niche (solving problem one), we may fall into the trap of comparing a point to the wrong elite, leading to an over or undervaluation of a point. In order to deal with this issue we calculate $\P(x \in c|D)$, the posterior probability of a point belonging to niche $c$. In this work we assume we have one feature, therefore $G(x)=g(x)$. We fit a second GP using data $X_n$ and the observed features $\{g(x_1),...,g(x_n)\}$ to make predictions $\bar{g}(x)$ and $s'(x)$ using Equations \ref{GP1} and \ref{GP2}. Given a point $x$, the probability that its feature $g(x)$ lies within a region $b_u< g(x) <b_l$, defining a niche $c$, can be calculated as
\begin{equation}
\P(x \in c|D) = \Phi\left(\frac{\bar{g}(x)-b_u}{ s'(x)}\right)-\Phi\left(\frac{\bar{g}(x)-b_l}{ s'(x)}\right),
\end{equation}
where $b_u$ and $b_l$ refer to the upper and lower boundary of the niche and $\bar{g}(x)$ and $s'(x)$ are the posterior mean and standard deviation from the feature GP at point $x$. This formulation can be trivially extended to calculate the probability of a point belonging to a niche in multidimensional feature space, though the test problems in this paper consider only one feature.

Putting these two pieces together we propose a niche specific Expected Improvement which returns the expected improvement over an elite, weighted by the posterior probability that it is in competition with that elite.
\begin{equation} \label{EIN}
\P(x \in c|D)EI_c(x) \label{nicheprob}
\end{equation}

A simple approach would be to sequentially optimise niches with this algorithm. This would certainly benefit from the shared surrogate models, allowing function calls for any single niche to improve the predictive power for all niches. However, our final acquisition function goes beyond this, estimating the expected improvement to the ensemble problem as a whole:
\begin{equation}
EJIE(x) = \sum_{i=1}^{C}\P(x\in c_i)EI_{c_i}(x) \label{EJIE}.
\end{equation}

The Expected Joint Improvement of Elites (EJIE) acquisition function gives high valuations to points that are likely to provide large improvements to the objective value of one or more niches. The predefined niches force the algorithm to search in areas of the input domain that single-objective optimisation would ignore, naturally diversifying the solution set. The posterior models provided by BOP-Elites provide insightful illumination of the search space, offering the opportunity for prediction of behavioural quality and objective performance of any point in the input domain with quantifiable confidence bounds. 

\begin{figure}
\centering
\begin{minipage}{.5\textwidth}
  \centering
  \includegraphics[width=\linewidth]{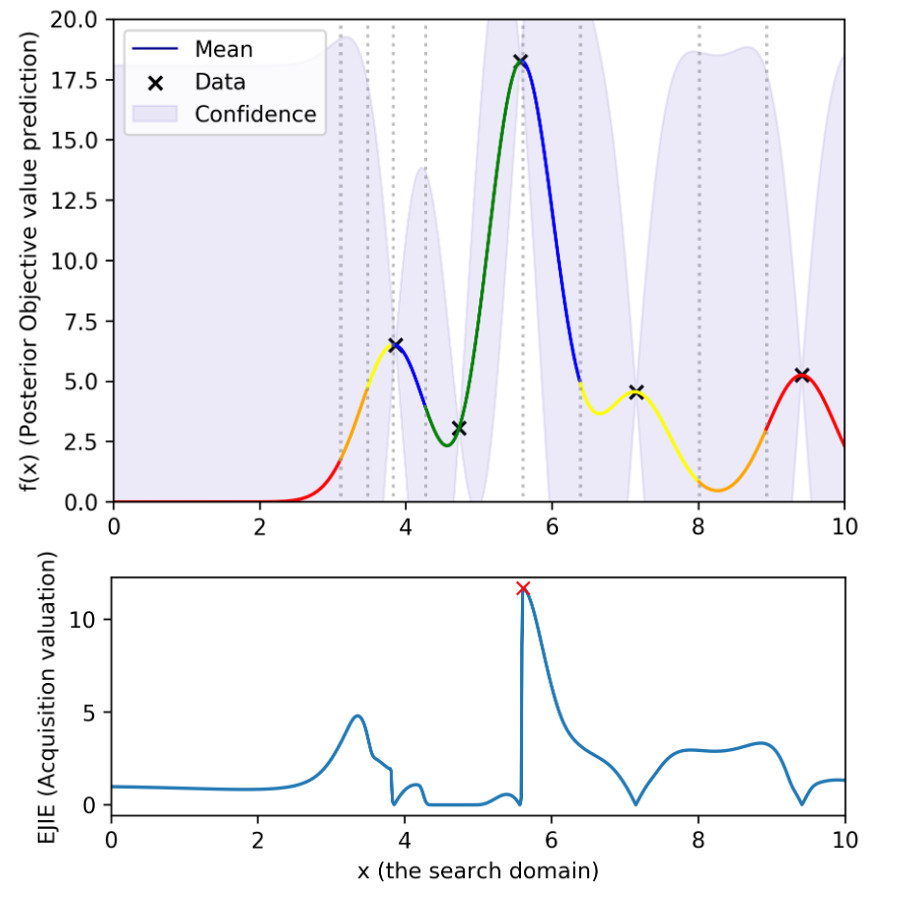}
  \label{fig:bopacq1}
\end{minipage}%
\begin{minipage}{.5\textwidth}
  \centering
  \includegraphics[width=\linewidth]{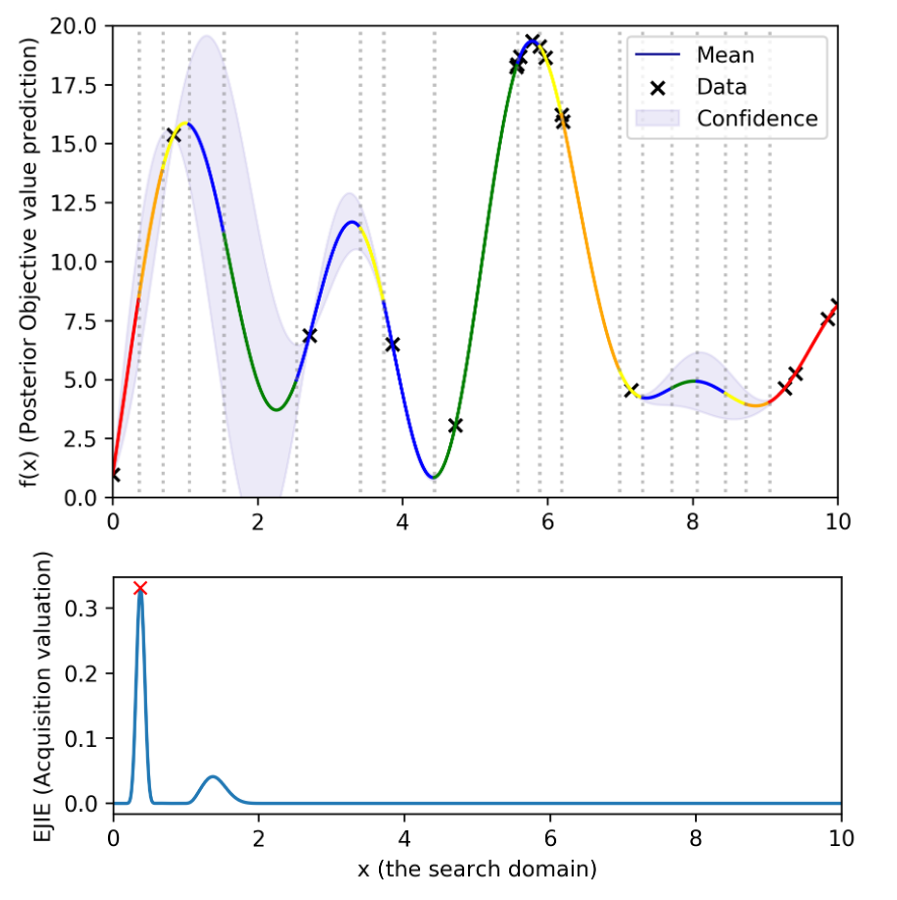}
  \label{fig:bopacq2}
\end{minipage}
\caption{A Plot of the posterior GP models showing the predicted objective value on the $y$ axis, the light blue area shows the uncertainty in the objective prediction. The posterior prediction of niche membership is indicated by the color of the line and class boundaries are indicated by the dotted vertical lines. Below each objective plot is the acquisition function plot, the suggested next point is highlighted with a red cross. [Left]: After five initial sample points there is a lot of uncertainty, EJIE predicts a large improvement at $x = 5.89$ [Right]: After 13 samples the model has approximated much of the objective function well and predicts a small improvement at $x = 0.24$ \label{posterioroutputfig} }
\end{figure}

\begin{algorithm}[H]
\caption{BOP-Elites Pseudo-Code}
\begin{algorithmic}[1]
 \Require $X$: Search domain 
 \Require $f(.)$: Objective function
 \Require $g(.)$: Feature function 
 \Require $N_{max}$: Sample budget 
 \Require $C$: List of Niches and their feature boundaries
 \Require $n_0$: Size of initial sampling budget
 \State Sample initial dataset $\{ X_{n_0},Y_{n_0},G(X_{n_{0}})\}=D_0$ randomly
 \State Create a set of observed elites, $\hat{S}$
 \State Create Objective and Feature surrogate model $\bar{f}(x),s(x),\bar{g}(x),s'(x)$ 
 \State  \textbf{For} $n=1$ \textbf{to} $N_{max}$ \textbf{do} 
 \State \indent \textbf{For} points in $X$ \textbf{do}: 
 \State \indent \indent compute all niche improvements (Eqn.~\ref{nicheimp}) $EI_{c_1}(x),...,EI_{c_C}(x)$
 \State \indent \indent compute all niche probabilities (Eqn.~\ref{nicheprob}) $\P(x\in c_1)...\P(x\in c_{C})$
 \State \indent \indent compute $EJIE(x)$ (Eqn.~\ref{EJIE})
 \State \indent sample $x_n$ at $argmax( EJIE(x))$
 \State \indent \textbf{if} $g(x_n) \in c_{i}$ and $f(x_n) > f(\hat{e}_{c_i})$ \textbf{do}
 \State \indent \indent $\hat{e}_{c_i} = x_n$
 \State \indent update posterior models $\bar{f}(x),s(x),\bar{g}(x),s'(x)$ \\
\Return{elites} 
\Statex{}
\end{algorithmic}
\end{algorithm}

\section{Empirical evaluation\label{sec:results}}
\subsection{Benchmark problem\label{problem_maker}}
100 benchmark problems consisting of an objective and feature function over a constrained, single dimensional search space $X = [0,10]$ were generated. Each problem was constructed in the following way: 11 random $y$ values were generated, uniformly at random from $Y = [0,20]$ and paired with the $x$ values [0,1,2...,10]. Independent GP's were fitted to these points and the posterior mean functions were taken. In this way we generated independent, random, problem functions of varying complexity dependent on the initial random points. 

It should be clear that it is possible for a feature function to be generated that does not pass through one or more niches. A badly designed algorithm would spend time looking for solutions that do not exist, and as we see later the EJIE acquisition function is able to avoid this particular issue.

\begin{figure}
\centering
\begin{minipage}{.5\textwidth}
  \centering
  \includegraphics[width=\linewidth]{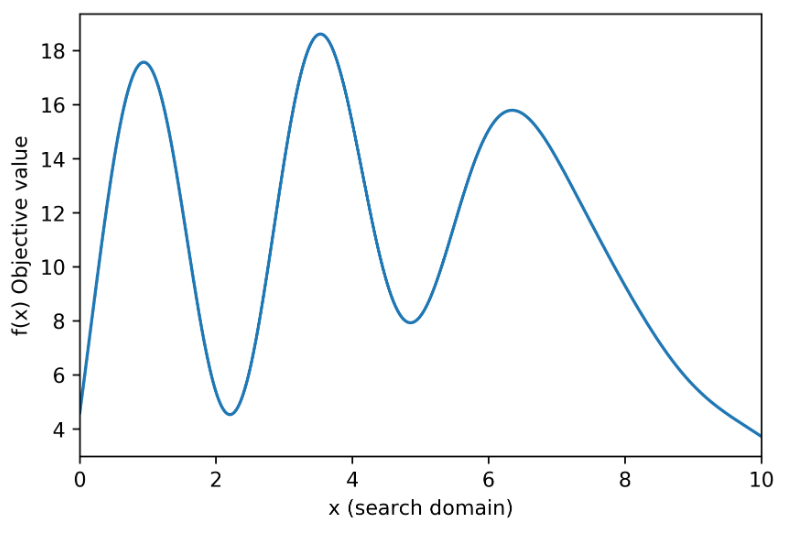}
  \label{fig:bpp1}
\end{minipage}%
\begin{minipage}{.5\textwidth}
  \centering
  \includegraphics[width=\linewidth]{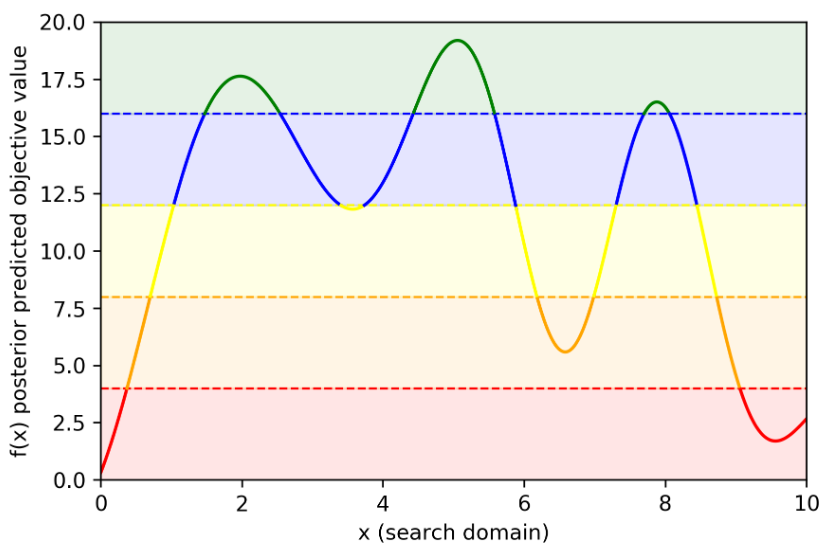}
  \label{fig:bpp2}
\end{minipage}
\caption{An example benchmark problem (same problem instance as used in Fig.~\ref{posterioroutputfig}), [Left]: A randomly generated objective function, [Right]: A random feature function with 5 niches }
\end{figure}

\subsection{Benchmark algorithm 1: Independent Niche Optimisation}
As a naive alternative to BOP-Elites we compare against an Independent Niche Optimisation algorithm following the same adaptation to EI as Eqn.~\ref{EIN}. This algorithm considers finding the optimal solution in each niche a constrained single-objective optimisation problem, splitting the search budget between the problems evenly and returning a single elite per niche. Instead of sharing surrogate models this algorithm wastefully builds independent models for each problem. We refer to this algorithm as IND in our results. 

\subsection{Benchmark algorithm 2: Sequential Niche Optimisation}
We consider the 'simple approach' briefly discussed in Section~\ref{acqfun}. This algorithm also uses Eqn.~\ref{EIN} as its acquisition function. The algorithm shares its search budget evenly between all niches and sequentially optimises each niche in turn (the ordering is considered arbitrary), searching for improvement over the elite for the current niche. Unlike the IND algorithm, all function calls are used to update shared surrogate models. In this sense, steps to improve the elite for any single niche improves the future predictive power of the surrogate models and therefore the optimisation of the entire problem. 

\subsection{Experimental setup}

We evaluated the performance of the 3 algorithms on 100 randomly generated test problems as described in Section~\ref{problem_maker} with $X=[0,10]$. The search domain was discretised into 1000 equally spaced points and 5 points were chosen uniformly at random. These points were evaluated on the objective and feature functions and the same 5 points given to each algorithm as their initial dataset. 

For all experiments objective and feature GPs were given an initial length-scale of 0.5, variance of 0.01 and a random restart budget of 100 for every iteration. Hyper-parameters were trained using maximum likelihood by L-BFGS. The length-scales were constrained to [0.001,2]. 

The metric used to evaluate solutions was TE, defined in Section~\ref{sec:probdef}. As we have ready access to the true functions we are able to easily compute the true optima required for this calculation by exhaustive search. In the following section, we plot this performance for every iteration of each algorithm.

\subsection{Experimental results}
\begin{figure}
\centering
\begin{minipage}{.5\textwidth}
  \centering
  \includegraphics[width=\linewidth]{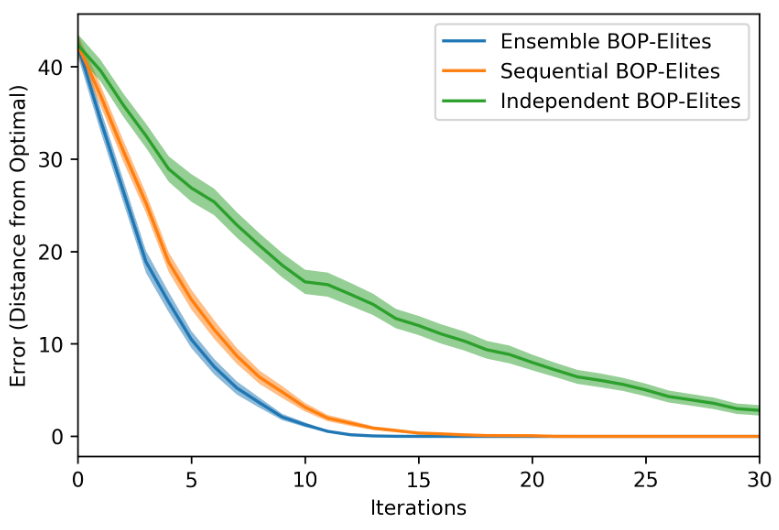}
  \label{fig:result1}
\end{minipage}%
\begin{minipage}{.5\textwidth}
  \centering
  \includegraphics[width=\linewidth]{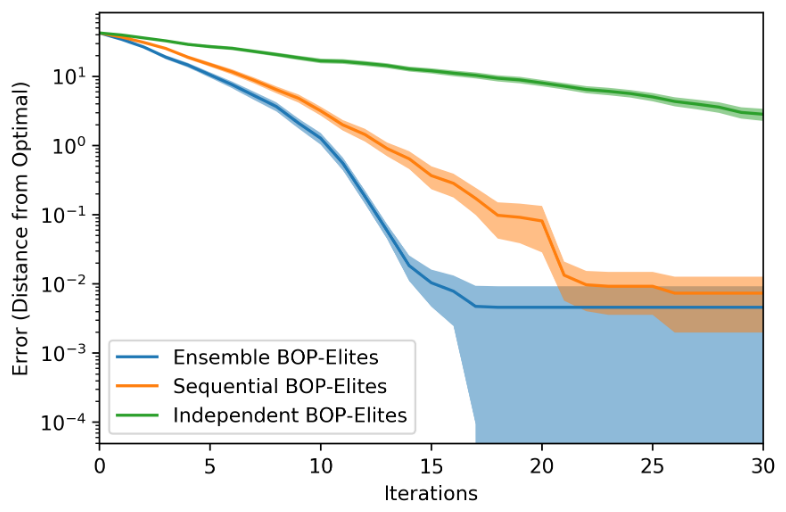}
  \label{fig:results2}
\end{minipage}
\caption{[Left]: A plot showing numerical results of the experiments for 3 algorithms over 100 random test functions. The Solid lines represent the mean performance, the shaded areas show the standard error around the mean [Right]: The same data plotted on a logscale to highlight the tail of the data \label{fig3}} 
\end{figure}

In our experiments, the ensemble BOP-Elites algorithm (blue lines in Fig.\ref{fig3}) converged, on average, to a near-optimal solution by the 13th iteration and the true optimal solution by the 18th iteration. BOP-Elites provided better overall solutions than both benchmark algorithms and drove towards the optimum faster.  

The only reason why BOP-Elites did not get a perfect score at the end is that in one out of the 100 replications BOP-Elites failed to find the optimal solution in one of the niches. In all other replications it reached the perfect solution within the available sampling budget.

The Sequential BOP-Elites algorithm (orange lines in Fig.~\ref{fig3}) performs well, but shows slower convergence to the optima. We note a visible 'kink' in the line of the log plot around iterations 17-20. This is likely explained by the sequential nature of the algorithm, causing the algorithm to keep sampling in niches for which the optimum has already been found, which leads to occasional slow-down. 

\section{Conclusion}
We have proposed a Bayesian optimisation algorithm for quality diversity search, in particular for identifying the best solution in each niche, where niche is defined through some additional feature function that needs to be explored.
The key idea is a novel acquisition function that uses Gaussian Process surrogate models for the objective as well as the feature function, and weighs the expected improvement with the probability of a solution belonging to each class.

The algorithm has demonstrated significantly faster convergence on a simple benchmark problem than alternative solution methods that allocate sampling budget to niches in a round robin fashion, or that consider the optimisation in each niche as a separate optimisation problem.

To the best of our knowledge, this is the first Bayesian optimisation algorithm for MAP-Elite problems. The algorithm can be further improved, e.g. by searching directly on the continuous space. It can be generalised to working on noisy problems and problems with higher dimensionality of search and feature spaces and should be tested on a wider range of problems, including more real-world problems.


\section*{ACKNOWLEDGEMENTS}
The authors would like to thank Jean-Baptiste Mouret for bringing this problem to our attention. The first author would like to acknowledge funding from EPSRC through grant EP/L015374/1.

\bibliographystyle{splncs04}
\bibliography{ref}




\end{document}